\documentclass[12pt, reqno]{amsart}
\usepackage{amsfonts,mathrsfs,amsmath,amssymb}

\def\1{\hbox{1\kern-.35em\hbox{1}}}

\newtheorem{theorem}{Theorem}[section]
\newtheorem*{theorem*}{Theorem}
\newtheorem{lem}[theorem]{Lemma}

\newtheorem*{proposition*}{Proposition}
\newtheorem{corollary}[theorem]{Corollary}

\newtheorem{remark}[theorem]{Remark}

\numberwithin{equation}{section}

\newcommand{\bea}{\begin{eqnarray}}
\newcommand{\eea}{\end{eqnarray}}
\newcommand{\be}{\begin{eqnarray*}}
\newcommand{\ee}{\end{eqnarray*}}

\newcommand{\Z}{{\mathbb Z}}

\newcommand{\C}{{\mathbb C}}

\def\tsv{\widetilde{\mathfrak{sv}}}

\newcommand{\Hom}{{\rm Hom}}

\def\qed{\hfill\mbox{$\Box$}}
\def\epsi{\epsilon}
\def\a{\alpha}

\def\d{\delta}

\def\al{\alpha}
\def\l{\lambda}

\def\si{\sigma}

\def\dis{\displaystyle}

\def\Z{\mathbb{Z}}

\def\C{\mathbb{C}}

\def\vp{\varphi}

\numberwithin{equation}{section}
\begin{document}
%
\title[derivations, central extensions and automorphism group of $W$] { The derivations, central extensions and
automorphism group of the Lie algebra $W$$^{*}$}
\author[Gao]{Shoulan Gao}
\address{Department of Mathematics, Shanghai Jiaotong University, Shanghai
200240, China} \email{gaoshoulan@sjtu.edu.cn}

\author[Jiang]{Cuipo Jiang$^{\dag}$}
\address{Department of Mathematics, Shanghai Jiaotong University, Shanghai
200240, China} \email{cpjiang@sjtu.edu.cn}

\author[Pei]{Yufeng Pei}
\address{Department of Mathematics, Shanghai Jiaotong University, Shanghai
200240, China} \email{yfpei@sjtu.edu.cn}

\begin{abstract}
In this paper,   we study the derivations, central extensions and
the automorphisms of the infinite-dimensional Lie algebra $W$ which
appeared in \cite{HSSU1} and Dong-Zhang's recent work \cite{ZD} on
the classification of some simple vertex operator algebras.

\end{abstract}

\thanks{$^{*}$ Supported in part by NSFC grant
10571119 and NSF grant 06ZR14049 of Shanghai City.}

\thanks{$^{\dag}$ Corresponding author: cpjiang@sjtu.edu.cn}
\maketitle

%
\section{\bf Introduction }
It is well known that the Virasoro algebra Vir plays an important
role in many areas of mathematics and physics (see \cite{KR1}, for
example). It can be regarded as the universal central extension of
the complexification of the Lie algebra Vect($S^1$) of (real) vector
fields on the circle $S^1$:
\begin{equation}
[L_m,L_n]=(m-n)L_{m+n}+\delta_{m+n,0}\frac{m^3-m}{12}c,
\end{equation}
where $c$ is a central element such that $[L_n,c] = 0$. The Virasoro
algebra admits many interesting extensions and generalizations, for
example, the $W_N$-algebras \cite{Z}, $W_{1+\infty}$ \cite{KR2}, the
higher rank Virasoro algebra \cite{PZ}, and the twisted
Heisenberg-Virasoro algebra \cite{ADKP,B1,JJ} etc. Recently M. Henke
et al. \cite{HSSU1,HSSU2}
 investigated a Lie algebra $W$ in their study of ageing phenomena which occur widely in physics \cite{H}.
 The Lie algebra $W$ is an abelian extension of centerless
Virasoro algebra, and is isomorphic to the semi-direct product Lie
algebra $\mathcal{L}\ltimes \mathcal{I}$, where $\mathcal{L}$ is the
centerless Virasoro algebra(Witt algebra) and $\mathcal{I}$ is the
adjoint $\mathcal{L}$-module. In particular, $W$ is the
infinite-dimensional extensions of the Poincare algebra
$\mathfrak{p}_3$. In their classification of the simple vertex
operator algebras with 2 generators, Dong and Zhang \cite{ZD}
studied a similar infinite-dimensional Lie algebra $W(2,2)$ and its
representation theory.  Although the algebra $W(2,2)$ is an
extension of the Virasoro algebra, as they remarked, the
representation theory for $W(2, 2)$ is totally different from that
of the Virasoro algebra.

The purpose of this paper is to study the structure of the Lie
algebra $W$ and its  universal central extension $\widetilde{W}$. We
will see that there is a natural surjective homomorphism from the
universal covering algebra $\widetilde{W}$ to $W(2,2)$.  We  show
that the second cohomology group with trivial coefficients for the
Lie algebra $W$  is two dimensional.
 Furthermore, by the Hochschild-Serre spectral sequence \cite{W} and R. Farnsteiner's theorem \cite{F},
we  determine the derivation algebras of $W$ and $\widetilde{W}$,
which  both have only one outer derivation. Finally, the
automorphism groups of $W$ and $\widetilde{W}$ are also
characterized.

Throughout the paper, we denote  by $\Z$  the set of all integers
and $\C$ the field of complex numbers.

%
\label{sub0-c-related}

%
\section{\bf The Universal Central Extension  of ${W}$}
%
\label{sub3-c-related}

The Lie algebra $W$ over the complex field $\C$ has a basis
$\{L_{m}, I_{m} \ | \ m\in\Z\}$ with the following bracket
\begin{eqnarray*}
&&[L_{m}, L_{n}]=(m-n)L_{m+n}, \ \ [L_{m}, I_{n}]=(m-n)I_{m+n},  \ \
[I_{m}, I_{n}]=0,
\end{eqnarray*}
for all $m,n\in\Z$.  It is clear that $W$ is isomorphic to the semi-direct
product Lie algebra $
W\simeq \mathcal{L}\ltimes\mathcal{I}$, where $\mathcal{L}=\bigoplus_{n\in\Z}\C L_{i}$ is
the classical  Witt algebra and $\mathcal{I}=\bigoplus_{n\in\Z}\C I_n$ can be regarded as
the adjoint $\mathcal{L}$-module.
Moreover, $W=\bigoplus\limits_{m\in\Z}W_{m}$ is a $\Z$-graded Lie algebra, where
$W_{m}=\C L_{m}\oplus\C I_{m}$.


Let $\mathfrak{g}$ be a  Lie algebra. Recall that a bilinear
function $\psi: \mathfrak{g}\times \mathfrak{g}\longrightarrow \C$
is called a 2-cocycle on $\mathfrak{g}$ if for all $x, y, z\in
\mathfrak{g}$, the following two conditions are satisfied:
\begin{equation}\nonumber
\psi(x, y)=-\psi(y, x),
\end{equation}
\begin{equation}\label{e2.0}
\psi([x, y], z)+\psi([y, z], x)+\psi([z, x], y)=0.
\end{equation}
 For any linear function $f:
\mathfrak{g}\longrightarrow \C$, one can define  a 2-cocycle
$\psi_{f}$ as follows
$$\psi_{f}(x, y)=f([x, y]), \quad\quad\quad \forall\; x, y\in \mathfrak{g}.$$
Such a 2-cocycle is called a 2-coboudary on $\mathfrak{g}$.

Let  $\mathfrak{g}$ be a  perfect Lie algebra, i.e., $[\mathfrak{g},
\mathfrak{g}]=\mathfrak{g}$. Denote by $C^{2}(\mathfrak{g}, \C)$ the
vector space of 2-cocycles on $\mathfrak{g}$, $B^{2}(\mathfrak{g},
\C)$ the vector space of 2-coboundaries on $\mathfrak{g}$. The
quotient space:
$$H^{2}(\mathfrak{g},\C)=C^{2}(\mathfrak{g},\C)/B^{2}(\mathfrak{g},\C)$$
is called the second cohomology group of $\mathfrak{g}$ with trivial
coefficients $\C$. It is well-known that $H^{2}(\mathfrak{g},\C)$
 is one-to-one correspondence to the equivalence classes of one-dimensional central extensions of the Lie
algebra $\mathfrak{g}$.

We will determine  the second cohomology group for the Lie algebra
${W}$.
\begin{lem}[See also \cite{OR}]Let $(\mathfrak{g},[\,,\,]_0)$ be a perfect Lie algebra over $\C$ and $V$
a $\mathfrak{g}$-module such that $\mathfrak{g}\cdot V=V$. Consider
the semi-direct product Lie algebra $(\mathfrak{g}\ltimes
V$,$[\,,\,])$ with the following bracket
$$[x,y]=[x,y]_0,\ \  [x,v]=x\cdot v,\ \  [u,v]=0, \quad\quad\forall\;x,y\in{\mathfrak{g}},\ \  u,v\in V.$$
Let $V^*$ be the dual $\mathfrak{g}$-module and
$$
B^{\mathfrak{g}}(V)=\{f\in \Hom(V\otimes V,\C)\,|\, f(u,v)=-f(v,u),\
\ f(x\cdot u,v)+f(u,x\cdot v)=0,\ \forall x\in\mathfrak{g},u,v\in
V\}.
$$
Then we have
$$H^2(\mathfrak{g}\ltimes V,\C)=H^2(\mathfrak{g},\C)\oplus
H^1(\mathfrak{g},V^*)\oplus B^{\mathfrak{g}}(V).$$
\end{lem}

\begin{proof}
Let $\al$ be a $2$-cocycle on $\mathfrak{g}\ltimes V$. Obviously,
$\al_{|\mathfrak{g}}\in H^2(\mathfrak{g},\C).$  Define $D_{\al}\in
\Hom_{\C}(\mathfrak{g}, V^*)$  by $$D_{\al}(x)(v)=\al(x,v),$$
 for
$x\in\mathfrak{g}, v\in V.$  Then   $D_{\al}\in
H^1(\mathfrak{g},V^*)$. In fact, for any
$x,y\in\mathfrak{g},v\in\;V$,
\begin{eqnarray*}
\al([x,y],v)&=&\al(x\cdot v ,y)-\al(y\cdot v,x)\\
&=&-D_{\al}(y)(x\cdot v)+D_{\al}(x)(y\cdot v)\\
&=&(x\cdot D_{\al}(y))(v)-(y\cdot D_{\al}(x))(v).
\end{eqnarray*}
Therefore
$$D_{\al}([x,y])=x\cdot D_{\al}(y)-y\cdot D_{\al}(x),$$
for all $x,y\in \mathfrak{g}$. Define $f_{\al}(u,v)=\al(u,v)$ for
any $u,v\in V$, then $f_{\al}\in B^{\mathfrak{g}}(V)$. It is straightforward to check  $\al_{|\mathfrak{g}}, D_{\al}$
and $f_{\al}$ are linearly independent. We get the desired formula.
\end{proof}


\begin{corollary} Let $V=\mathfrak{g}$ be  the adjoint $\mathfrak{g}$-module. Then
\begin{eqnarray}
H^2(\mathfrak{g}\ltimes
\mathfrak{g},\C)&=&H^2(\mathfrak{g},\C)\oplus
H^1(\mathfrak{g},\mathfrak{g}^*)\oplus B^{\mathfrak{g}}(\mathfrak{g}).
\end{eqnarray}
\end{corollary}

\begin{theorem}\label{extension}
$H^2(W,\C)=\C\alpha\oplus\C\beta,$ where
\begin{eqnarray*}
&&\alpha(L_m,L_n)=\delta_{m+n,0}\frac{m^3-m}{12},\quad \alpha(L_m,I_n)=\alpha(I_m,I_n)=0,\\
&&\beta(L_m,I_n)=\delta_{m+n,0}\frac{m^3-m}{12},\quad
\beta(L_m,L_n)=\beta(I_m,I_n)=0.
\end{eqnarray*}
for any $m,n\in\Z$.
\end{theorem}
\begin{proof}
Since $H^2(\mathcal{L},\C)=\C\alpha$ and $
H^1(\mathcal{L},\mathcal{L}^*)\simeq \C\beta$ (see \cite{LP, HPL}
for the details), we need to prove that
$B^{\mathcal{L}}(\mathcal{I})=0$. Let $f\in
B^{\mathcal{L}}(\mathcal{I})$, then
\begin{equation}\label{2.1}
(i-j)f(I_{i+j}, I_{k})+(k-i)f(I_{k+i},  I_{j})=0.
\end{equation}
Letting $i=0$ in (\ref{2.1}), we get
$$(j+k)f(I_{j}, I_{k})=0.$$
So $f(I_{j}, I_{k})=0$ for $j+k\neq 0$. Let $k=-i-j$ in (\ref{2.1}),
then we obtain
\begin{equation}\label{2.2}
(i-j)f(I_{i+j}, I_{-i-j})+(2i+j)f(I_{j},  I_{-j})=0.
\end{equation}
Let $j=-i$, then $2if(I_{0}, I_{0})=if(I_{i}, I_{-i})$,  which
implies that $f(I_{i}, I_{-i})=0$ for all $i\in\Z$. Therefore,
$f(I_{m}, I_{n})=0$ for all $m, n\in\Z$.

\end{proof}


Let $\mathfrak{g}$ be a Lie algebra,  $(\widetilde{\mathfrak{g}},
\pi)$ is called a central extension of $\mathfrak{g}$ if $\pi:
\widetilde{\mathfrak{g}}\longrightarrow \mathfrak{g}$ is a
surjective homomorphism whose kernel lies in the center of the Lie
algebra $\widetilde{\mathfrak{g}}$. The pair
$(\widetilde{\mathfrak{g}}, \pi)$ is called a covering of
$\mathfrak{g}$ if  $\widetilde{\mathfrak{g}}$ is perfect. A covering
$(\widetilde{\mathfrak{g}}, \pi)$ is called a universal central
extension of $\mathfrak{g}$ if for every central extension
$(\widetilde{\mathfrak{g}}', \varphi)$ of $\mathfrak{g}$ there is a
unique homomorphism $\psi:
\widetilde{\mathfrak{g}}\longrightarrow\widetilde{\mathfrak{g}}'$
for which $\varphi\psi=\pi$. It follows from \cite{H} that every
perfect Lie algebra has a universal central extension.

Let $\widetilde{{W}}={W}\bigoplus\C\;C_{1}\bigoplus\C\;C_{2}$ be a
vector space over the complex field $\C$ with a basis $\{L_{n},
I_{n}, C_{1}, C_{2} \ |\ n\in\Z \}$ satisfying the following
relations
$$[L_{m}, L_{n}]=(m-n)L_{m+n}+\delta_{m+n,0}\dis\frac{m^{3}-m}{12}C_{1},$$
$$[L_{m}, I_{n}]=(m-n)I_{m+n}+\delta_{m+n,0}\dis\frac{m^{3}-m}{12}C_{2},$$
$$[I_{m}, I_{n}]=0,\ \ [C_{1}, \widetilde{{W}}]=0,\ \ [C_{2}, \widetilde{{W}}]=0, $$
for all $m,n\in\Z$. By Theorem \ref{extension}, $\widetilde{{W}}$ is
a universal covering algebra of ${W}$.

There is  a $\Z$-grading on $\widetilde{{W}}$:
$$\widetilde{{W}}=\bigoplus\limits_{n\in\Z}\widetilde{{W}}_{n},$$
where $\widetilde{W}_{n}=span\{L_{n}, I_{n}, \d_{n, 0}C_{1}, \d_{n,
0}C_{2} \}$. Set
$\widetilde{W}_{+}=\bigoplus\limits_{n>0}\widetilde{{W}}_{n}$ and
$\widetilde{W}_{-}=\bigoplus\limits_{n<0}\widetilde{{W}}_{n}$,
then there is a triangular decomposition on $\widetilde{{W}}$:
$$\widetilde{{W}}=\widetilde{W}_{-}\bigoplus\widetilde{W}_{0}\bigoplus\widetilde{W}_{+}.$$

\begin{remark}
Let $C_1=C_2=C$, then we get the Lie algebra $W(2,2)$ appeared in
\cite{ZD}.
\end{remark}

%
\section{\bf The Derivation Algebra of $W$}
%
\label{sub3-c-related}

In this section, we consider the derivation algebra of $W$.

Let $G$ be a commutative group, $\mathfrak{g}=\bigoplus\limits_{g\in
G}\mathfrak{g}_{g}$ a $G$-graded Lie algebra. A
$\mathfrak{g}$-module $V$ is called $G$-graded, if
$$V=\bigoplus\limits_{g\in G}V_{g\in G},
\;\;\; \mathfrak{g}_{g}V_{h}\subseteq V_{g+h},\;\;\; \forall\;
g,h\in G.$$ Let $\mathfrak{g}$ be a Lie algebra and $V$  a
$\mathfrak{g}$-module. A linear map $D: \mathfrak{g}\longrightarrow
V$ is called a derivation, if for any $x, y\in \mathfrak{g}$,
$$D[x, y]=x.D(y)-y.D(x).$$
If there exists some $v\in V$ such that $D:x\mapsto x.v$, then $D$
is called an inner derivation. Denote by $Der(\mathfrak{g}, V)$ the
vector space of all derivations, $Inn(\mathfrak{g}, V)$ the vector
space of all inner derivations. Set
$$H^{1}(\mathfrak{g}, V) =Der(\mathfrak{g}, V) /Inn(\mathfrak{g}, V).$$
Denote by $Der(\mathfrak{g})$ the derivation algebra of
$\mathfrak{g}$, $Inn(\mathfrak{g})$  the vector space of all inner
derivations of $\mathfrak{g}$.

Firstly, we have the short exact sequence of Lie algebras
$$
0\to \mathcal{I}\to W\to \mathcal{L}\to 0,
$$
which induces an exact
sequence
\begin{equation}
H^1(W,\mathcal{I})\to H^1(W,W)\to H^1(W,\mathcal{L}).
\end{equation}
The right-hand side can be computed from the initial terms
\begin{equation}
0\to H^1(\mathcal{L},\mathcal{L})\to H^1(W,\mathcal{L})\to H^1(\mathcal{I},\mathcal{L})^W
\end{equation}
of the four-term sequence associated to the Hochschild-Serre
spectral sequence \cite{W}. It is known that
$H^1(\mathcal{L},\mathcal{L})=0$, while the term
$H^1(\mathcal{I},\mathcal{L})^W$ is equal to
$\Hom_{U(W)}(\mathcal{I}/[\mathcal{I},\mathcal{I}],\mathcal{L})$.
Therefore, it suffices to compute the terms $H^1(W,\mathcal{I})$ and
$\Hom_{U(W)}(\mathcal{I}/[\mathcal{I},\mathcal{I}],\mathcal{L})$.

By Proposition 1.1 in \cite{F}, we have the following lemma.
\begin{lem}\label{L3.2}
$$Der(W,\mathcal{I})=\bigoplus\limits_{n\in\Z}Der(W,\mathcal{I})_{n},$$
where $Der(W,\mathcal{I})_{n}(\mathcal{I}_{m})\subseteq \mathcal{I}_{m+n}$ for all $m,n\in\Z$.
\end{lem}

\qed

\begin{lem}\label{L3.3} For any nonzero integer $m$, we have
$$H^{1}(W_{0},\mathcal{I}_{m})=0,$$ where $\mathcal{I}_{m}=\C I_{m}$.
\end{lem}

\begin{proof}  Let $m\neq 0$,  $\vp:W_{0}\longrightarrow
\C I_{m} $  a derivation. Assume that
\begin{eqnarray*}
& &\vp(L_{0})=aI_{m}, \ \ \vp(I_{0})=bI_{m},
\end{eqnarray*}
where $a,b\in\C$. Since $\vp[L_{0},I_{0}]=[\vp(L_{0}),I_{0}]+[L_{0},
\vp(I_{0})],$ it is easy to deduce $b=0$. Let
$E_{m}=-\dis\frac{a}{m}I_{m}$,  then we have $\vp(L_{0})=[L_{0},
E_{m}]$ and $\vp(I_{0})=[I_{0}, E_{m}].$ Therefore, $\vp\in
Inn(W_{0},I_{m})$, that is,
$$H^{1}(W_{0},I_{m})=0, \quad \forall\;
m\in\Z\setminus\{0\}.$$
\end{proof}

\begin{lem}\label{L3.4}
$\Hom_{ W_{0}}( W_{m},  \mathcal{I}_n)=0$ for all $m,n\in\Z$, $m\neq
n.$
\end{lem}

\begin{proof} Let $f\in \Hom_{ W_{0}}(W_{m}, \mathcal{I}_n)$, where $m\neq n$. Then
for any $ E_{m}\in W_{m}$, we have
$$f([L_{0}, E_{m}])=[L_{0}, f(E_{m})].$$ Note $[L_{0},
E_{k}]=-kE_{k}$ for all $k\in\Z$. Then  we get
$$-mf( E_{m})=[L_{0}, f(E_{m})]=-nf( E_{m}).$$
So  $f( E_{m})=0$ for all $m\neq n$. Therefore, we have $f=0$.
\end{proof}

By Lemma \ref{L3.3}-\ref{L3.4} and Proposition 1.2 in \cite{F}, we
have the following Lemma.

\begin{lem}\label{L3.5}
$Der(W,\mathcal{I})=Der(W,\mathcal{I})_{0}+Inn(W,\mathcal{I}).$
\end{lem}

\begin{lem}\label{L3.6}
$H^{1}(W,\mathcal{I})=\C D_{1}$, where
$$D_{1}(L_{m})=0,\ \   D_{1}(I_{m})=I_{m},$$ for all $m\in\Z$ .
\end{lem}
\begin{proof}

For any $D\in Der(W,\mathcal{I})_{0}$, assume
$$D(L_{m})=a_{m}I_{m}, \ \  D(I_{m})=b_{m}I_{m},$$
where $a_{m}, b_{m}\in\C$. By the definition of derivation and the
Lie bracket  in $W$, we have
$$a_{m+n}=a_{m}+a_{n},  \ \   b_{m+n}=b_{n},  \ \ \ \ m\neq n.$$
Obviously, $b_{m}=b_{0}$ for all $m\in\Z$ and $a_{0}=0,
a_{-m}=-a_{m}$. By induction on $m>0$, we  deduce that
$a_{m}=(m-2)a_{1}+a_{2}$ for $m>2$.  Then it is easy to infer that
$a_{2}=2a_{1}$. Consequently, we get $a_{m}=ma_{1}$ for all
$m\in\Z$. So for all $m\in\Z$, we have
$$D(L_{m})=ma_{1}I_{m}, \ \  D(I_{m})=b_{1}I_{m}.$$
 Set
$D_{0}=-ad(a_{1}I_{0})\in{Inn(W)}$, then $D(L_{m})=D_{0}(L_{m}),\ \
D(I_{m})=D_{0}(I_{m})+b_{1}I_{m}.$ Therefore, $$\bar{D}(L_{m})=0,
\bar{D}(I_{m})=b_{1}I_{m},$$ for all $m\in\Z$, where
$\bar{D}=D-D_{0}$ is an outer derivation. The lemma holds.

\end{proof}

\begin{lem}\label{L3.6}
$\Hom_{U(W)}(\mathcal{I}/[\mathcal{I},\mathcal{I}],\mathcal{L})=0.$
\end{lem}
\begin{proof}
 As a matter of fact,
 $$\Hom_{U(W)}(\mathcal{I}/[\mathcal{I},\mathcal{I}],\mathcal{L})=\Hom_{U(W)}(\mathcal{I},\mathcal{L}).$$
For any $f\in \Hom_{U(W)}(\mathcal{I},\mathcal{L})$, we have
$[L_{0}, f(I_{m})]=f([L_{0}, I_{m}])$, i.e.,
$$ad(-L_{0})(f(I_{m}))=mf(I_{m}),$$
for all $m\in\Z$. This suggests $f(I_{m})\in \C L_{m}$. Assume
$f(I_{m})=x_{m}L_{m}$ for all $m\in\Z$, where $x_{m}\in\C$. By the
relation  that $[L_{n}, f(I_{m})]=f([L_{n}, I_{m}])$ for all
$m,n\in\Z$, we have
$$x_{m+n}=x_{m},  \ \ \ \ m\neq n.$$
Obviously, $x_{m}=x_{0}$ for all $m\in\Z$. So there exists some
constant $a\in\C$ such that $$f(I_{m})=aL_{m},$$ for all $m\in\Z$.
Since  $f([I_0,I_1])=0=[I_0,f(I_1)]=-aL_1,$ we have $a=0$. Hence
$f=0$.

\end{proof}

\begin{theorem}
$H^{1}(W,W)=\C D$, where
$$D(L_{m})=0,\ \   D(I_{m})=I_{m},$$ for all $m\in\Z$ .
\end{theorem}

Furthermore, it follows from Theorem 2.2 in \cite{BM} that
\begin{corollary}
$H^1(\widetilde{{W}},\widetilde{{W}})\simeq H^{1}(W,W).$
\end{corollary}

%
\section{\bf  The Automorphism Group of  ${W}$ }
%
\label{sub5-c-related}

 Denote by $Aut({W})$ and $\mathfrak{Z}$ the
automorphism group and the inner automorphism group of ${W}$
respectively. Obviously, $\mathfrak{Z}$  is generated by $\exp(k ad
I_{m})$, $m\in\Z,$ $k\in\C$, and $\mathfrak{Z}$ is an abelian
subgroup.
%
Note that $\mathcal{I}$  is the maximal proper ideal of $W$, so we
have the following lemma.

\begin{lem}\label{L3.1}
For any $\si\in Aut({W})$, $\si(I_{n})\in \mathcal{I}$ for all
$n\in\Z$.
\end{lem}

\qed

For any $\prod\limits_{j=s}^{t}\exp(k_{i_{j}}adI_{i_{j}})\in
\mathfrak{Z}$, we have
$$\prod\limits_{j=s}^{t}\exp(k_{i_{j}}adI_{i_{j}})(I_{n})=I_{n},\;\;
\prod\limits_{j=s}^{t}\exp(k_{i_{j}}adI_{i_{j}})(L_{n})=L_{n}+\sum\limits_{j=s}^{t}k_{i_{j}}(i_{j}-n)I_{i_{j}+n}.$$

\begin{lem}\label{L3.2}
For any $\si\in Aut({W})$, there exist some $\tau\in \mathfrak{Z}$
and $\epsi\in\{\pm1\}$ such that
\begin{equation}\label{E3.1}
\bar{\si}(L_{n})= a^{n}\epsi L_{\epsi n}+a^{n}\lambda n I_{\epsi n},
\end{equation}
\begin{equation}\label{E3.2}
\bar{\si}(I_{n})= a^{n}\mu I_{\epsi n},
\end{equation}
where $\bar{\si}=\tau^{-1}\si$, $a, \mu\in\C^{*}$ and $\l\in\C$.
Conversely, if $\bar{\si}$ is a linear operator on ${\tsv}$
satisfying (\ref{E3.1})-(\ref{E3.2}) for some $\epsi\in\{\pm1\}$,
$a,\mu\in\C^{*}$ and $\l\in\C$, then $\bar{\si}\in Aut({W})$.
\end{lem}

\begin{proof} For any $\si\in Aut({W})$, denote
$\si|_{\mathcal{L}}=\si'$. Then  $\si'$ is an automorphism of the
classical Witt algebra, so $\si'(L_{m})=\epsilon a^{m}L_{\epsilon
m}$ for all $m\in\Z$, where $a\in\C^{*}$ and $\epsi\in\{\pm1\}$.
Assume that
$$\si(L_{0})= \epsi L_{0}+\sum\limits_{i=p }^{q} \lambda_{i} I_{i}+\lambda_{0} I_{0},$$
where $i\neq 0$. Let
$\tau=\prod\limits_{i=p}^{q}\exp(\frac{\l_{i}}{\epsi i}adI_{i})$,
then
$$\tau(\epsi L_{0})=\epsi L_{0}+\sum\limits_{i=p }^{q} \lambda_{i} I_{i}.$$
Therefore, $\si(L_{0})= \tau(\epsi L_{0})+\lambda_{0} I_{0}.$ Set
$\bar{\si}=\tau^{-1}\si$, then $\bar{\si}(L_{0})=\epsi
L_{0}+\lambda_{0} I_{0}.$ Assume
$$\bar{\si}(L_{n})= a^{n}\epsi L_{\epsi
n}+a^{n}\sum\lambda(n_{i})I_{n_{i}}, \;\; n\neq 0,$$ where each
formula is of finite terms and $ \lambda(n_{i}), \mu(n_{j})\in\C$.
For $m\neq 0$, through  the relation
$[\bar{\si}(L_{0}),\bar{\si}(L_{m})]=-m\bar{\si}(L_{m})$, we get
$$\l_{0}mI_{\epsi m}=\sum\l(m_{i})[m-\epsi m_{i}]I_{m_{i}}.$$
This forces that $\l_{0}=0.$ Then
$[\bar{\si}(L_{0}),\bar{\si}(L_{m})]=[\epsi
L_{0},\bar{\si}(L_{m})]=-m\bar{\si}(L_{m})$, that is,
$$-adL_{0}(\bar{\si}(L_{m}))=\epsi m\bar{\si}(L_{m}).$$ So for all
$m\in\Z$,  we obtain $$\bar{\si}(L_{m})= a^{m}\epsi L_{\epsi
m}+a^{m}\lambda(\epsi m)I_{\epsi m}.$$ Comparing the coefficients of
$I_{\epsi(m+n)}$ on the both side of
$[\bar{\si}(L_{m}),\bar{\si}(L_{n})]=(m-n)\bar{\si}(L_{m+n})$, we
get $$\l(\epsi m)+\l(\epsi n)=\l(\epsi (m+n)).$$ Note $\l_{0}=0$, we
 deduce that $\l(\epsi m)=m\l(\epsi )$ for all $m\in\Z$. Therefore,
$$\bar{\si}(L_{m})= a^{m}\epsi L_{\epsi\;m}+a^{m}m\lambda(\epsi )I_{\epsi m}.$$
Finally,  by  $[{\si}(L_{0}), {\si}(I_{m})]=-m{\si}(I_{m})$, we have
$$-adL_{0}(\si(I_{m}))=\epsi m \si(I_{m}),$$ for all $m\in\Z$. Then
by Lemma \ref{L3.1}, we may assume
$$\si(I_{n})= a^{n}\mu(\epsi n) I_{\epsi n},$$
for all $n\in\Z$. By $[{\si}(L_{m}),
{\si}(I_{n})]=(m-n){\si}(I_{m+n})$, we get
$$\mu(\epsi n)=\mu(\epsi (m+n)),\ \ \ \ m\neq n.$$
Consequently, $\mu(\epsi m)=\mu(0)$ for all $m\in\Z$. Set
$\mu(0)=\mu$, then for all $n\in\Z$, we have
$$\bar{\si}(I_{n})= a^{n}\mu I_{\epsi n}.$$
\end{proof}

Denote by $\bar{\si}(\epsi, \lambda, a, \mu)$ the automorphism of
$W$ satisfying (\ref{E3.1})-(\ref{E3.2}), then
\begin{equation}\label{5.30}
\bar{\si}(\epsi_{1}, \lambda_{1}, a_{1}, \mu_{1})
\bar{\si}(\epsi_{2},\lambda_{2},a_{2},\mu_{2})
=\bar{\si}(\epsi_{1}\epsi_{2}, \lambda_{1}+\mu_{1}\lambda_{2},
a_{1}^{\epsi_{2}}a_{2}, \mu_{1}\mu_{2}),
\end{equation}
and $\bar{\si}(\epsi_{1}, \lambda_{1}, a_{1},
\mu_{1})=\bar{\si}(\epsi_{2}, \lambda_{2}, a_{2}, \mu_{2})$ if and
only if $\epsi_{1}=\epsi_{2}, \lambda_{1}=\lambda_{2}, a_{1}=a_{2},
\mu_{1}= \mu_{2}$. Let
$$\bar{\pi}_{\epsi}=\bar{\si}(\epsi, 0, 1, 1),\quad
\bar{\si}_{\lambda}=\bar{\si}(1, \lambda, 1, 1),\quad \bar{\si}_{a,
\mu}=\bar{\si}(1, 0, a, \mu)$$ and
$$\mathfrak{a}=\{\bar{\pi}_{\epsi} \;|\; \epsi=\pm1 \},\quad
\mathfrak{t}=\{\bar{\si}_{\lambda} \;|\; \lambda\in\C \},\quad
\mathfrak{b}=\{\bar{\si}_{a,\mu} \; |\; a,\mu\in\C^{*} \}. $$ By
(\ref{5.30}), we have the following relations:
$$\bar{\si}(\epsi, \lambda, a, \mu)=\bar{\si}(\epsi, 0,
1,1)\bar{\si}(1, \lambda, 1, 1)\bar{\si}(1, 0, a, \mu)\in
\mathfrak{a}\mathfrak{t}\mathfrak{b}  ,$$
$$\bar{\si}(\epsi, \lambda, a, \mu)^{-1}=\bar{\si}(\epsi, -\lambda\mu^{-1}, a^{-\epsi},\mu^{-1})  ,$$
$$\bar{\pi}_{\epsi_{1}}\bar{\pi}_{\epsi_{2}}=\bar{\pi}_{\epsi_{1}\epsi_{2}},\quad
 \bar{\si}_{\lambda_{1}}\bar{\si}_{\lambda_{2}}=\bar{\si}_{\lambda_{1}+\lambda_{2}},\quad
\bar{\si}_{a_{1},\mu_{1}}\bar{\si}_{a_{2},\mu_{2}}=\bar{\si}_{a_{1}a_{2},\mu_{1}\mu_{2}},$$
$$\bar{\si}_{\lambda}\bar{\pi}_{\epsi}=\bar{\pi}_{\epsi}\bar{\si}_{\lambda}, \quad\;
\bar{\pi}_{\epsi}^{-1}\bar{\si}_{a,\mu}\bar{\pi}_{\epsi}=\bar{\si}_{a^{\epsi},\mu},\quad
\bar{\si}_{a,\mu}\bar{\si}_{\lambda}\bar{\si}_{a,\mu}^{-1}=\bar{\si}_{\mu\l}.$$
Hence,  the following lemma holds.
\begin{lem}\label{L3.3}
$\mathfrak{a},\mathfrak{t}$ and $ \mathfrak{b} $ are all subgroups
of $Aut({W})$.  Furthermore, $\mathfrak{t}$  is an abelian normal
subgroup  commutative with $\mathfrak{Z}$.
$$Aut({W})=({\mathfrak{Z}\mathfrak{t}})\rtimes(\mathfrak{a}\ltimes\mathfrak{b}),$$
where $\mathfrak{a}\cong\Z_{2}=\{\pm 1\}, \mathfrak{t}\cong\C $,
$\mathfrak{b}\cong\C^{*}\times\C^{*}.$ \qed
\end{lem}

Let $\C^{\infty}=\{(a_{i})_{i\in\Z}\;|\; a_{i}\in\C, {\rm{\; all \;
but\;  \; a \;  finite \; number \; of \;  the}} \; a_{i} \; {\rm{
\; are \; zero}} \ \}$. Then $\C^{\infty}$ is an abelian group.

\begin{lem}\label{L3.4}
 ${\mathfrak{Z}\mathfrak{t}}$ is isomorphic to $\C^{\infty}$.
\end{lem}
\begin{proof}
Define $f: {\mathfrak{Z}\mathfrak{t}}\longrightarrow \C^{\infty}$ by
$$f(\bar{\si}_{\lambda}\prod\limits_{i=1}^{s}exp(\a_{k_{i}}{\rm{ad}} I_{k_{i}}))=(a_{p})_{p\in\Z},$$
 where $a_{k_{i}}=\a_{k_{i}}$ for $k_{i}<0$,  $a_{0}=\l$,   $a_{k_{i}+1}=\a_{k_{i}}$ for $k_{i}\geq0$, and the others are zero, $k_{i}\in\Z$ and $k_{1}<
k_{2}<\cdots< k_{s}$. Since every element of ${{\rm
\mathfrak{Z}}\mathfrak{t}}$ has the unique form of
$\bar{\si}_{\lambda}\prod\limits_{i=1}^{s}exp(\a_{k_{i}}{\rm{ad}}
I_{k_{i}})$, it is easy to check that $f$ is an isomorphism of
group.

\end{proof}

\begin{theorem}\label{T3.5}
$Aut({W})\cong\C^{\infty}\rtimes(\Z_{2}\ltimes(\C^{*}\times\C^{*})).$
\end{theorem}
\qed

Since $W$ is centerless,  it follows from Corollary 6 in \cite{P}
that $Aut({\widetilde{W}})=Aut({W})$, that is,
$$Aut({\widetilde{W}})\cong\C^{\infty}\rtimes(\Z_{2}\ltimes(\C^{*}\times\C^{*})).$$

\bibliography{}

\end{document}